\documentclass[11pt]{amsart}

\usepackage {amsmath, amssymb, amsthm}
\usepackage[ruled,vlined]{algorithm2e}
\usepackage [top=1.5in, right=1in, bottom=1.5in, left=1in] {geometry}
\usepackage {pgfplots}
\pgfplotsset {compat=1.13}
\usepackage {tikz}
\usepackage{todonotes}
\usepackage{cite}
\usepackage{graphicx}
\usepackage {hyperref}
\usepackage{enumitem}
\newcommand{\mylabel}[2]{#2\def\@currentlabel{#2}\label{#1}}

\hypersetup {
    colorlinks = true,
    linktoc = all,
    linkcolor = blue,
    citecolor = red
}
\usepackage{xcolor}

\usepackage {libertine}
\usepackage [T1]{fontenc}

\theoremstyle {plain}
\newtheorem {thm} {Theorem}[section]
\newtheorem {prop}[thm] {Proposition}

\newtheorem {lem}[thm] {Lemma}

\theoremstyle {definition}

\theoremstyle {remark}
\newtheorem {rmk}[thm] {Remark}

\newcommand{\Z}{\mathbb{Z}}

\newcommand{\w}{\mathbf{w}}
\newcommand{\x}{\mathbf{x}}
\newcommand{\A}{\mathbf{a}}
\newcommand{\B}{\mathbf{b}}
\newcommand{\C}{\mathbf{c}}
\newcommand{\vtm}{\mathbf{vtm}}
\newcommand{\tm}{\mathbf{tm}}
\newcommand{\pd}{\mathbf{pd}}
\renewcommand{\u}{\mathbf{u}}
\renewcommand{\v}{\mathbf{v}}
\newcommand{\rev}{\mathrm{rev}}
\newcommand{\rot}{\mathrm{rot}}

\begin {document}

\title{Group Actions and Some Combinatorics on Words with $\vtm$}

\author{John Machacek}
\email{jmachacek@hsc.edu}
\address{Department of Mathematics and Computer Science, Hampden-Sydney College}

\begin{abstract}
We introduce generalizations of powers and factor complexity via orbits of group actions.
These generalizations include concepts like abelian powers and abelian complexity.
It is shown that this notion of factor complexity cannot be used to recognize Sturmian words in general.
Within our framework, we establish square avoidance results for the ternary squarefree Thue--Morse word $\vtm$.
These results go beyond the usual squarefreeness of $\vtm$ and are proved using Walnut.
Lastly, we establish a group action factor complexity formula for $\vtm$ that is expressed in terms of the abelian complexity of the period doubling word $\pd$.
\end{abstract}
    
\maketitle

\section{Introduction}

Combinatorics on words studies finite and infinite sequences of symbols.
The origins~\cite{origins} of the field are Thue's constructions on avoiding repetitions in words~\cite{berstel1992axel}.
One of the words considered by Thue is the word that we will denote $\vtm$ (also called the Hall--Thue word or the ternary squarefree Thue--Morse word), which will serve as a central example in this paper.
A key development in combinatorics on words was the work Morse--Hedlund~\cite{MorseHedlund} on symbolic dynamics that, among other things, developed the notion of counting distinct factors (i.e., contiguous subsequences) of a word as a measure of complexity.
These classic concepts of repetition and factor complexity are what we aim to study from a more generalized perspective.

Words are ordered structures by their nature.
That is, the symbols do not commute by default.
When one allows symbols to commute, we get the subfield of abelian combinatorics on words (see, e.g.,~\cite{abelian}).
Letting symbols commute can be thought of as letting orbits of a word define equivalence classes, where here we have an action of a symmetric group that permutes positions.
One natural direction of generalization is to use subgroups of the symmetric group to permute positions rather than full symmetric groups.
This is done in~\cite{GroupComplexity}, where such group actions and counting equivalence classes of factors were considered to prove variations of the Morse--Hedlund theorem~\cite{MH40}.
Our setting is similar, with a difference being that we also allow an additional group acting by permuting our alphabet.

We will use equivalence classes defined by orbits to look at generalizations of repetitions in words.
One of the most studied types of repetition is squares, which is a word that is the concatenation of a shorter word with itself.
For example, the English words \emph{meme} or \emph{couscous} are squares.
The word $\vtm$ famously does not contain any squares.
A \emph{square} is a word of the form $\x^2 = \x\x$, and we will study words of the form $\x\x'$ where $\x$ and $\x'$ are in the same orbit under some group action.
Examples of the latter type of repetition are the English words \emph{reappear} or \emph{teammate}.
Both of these examples are what are known as abelian squares because the second half of each is a permutation of the first half.
We can refine our approach group theoretically and differentiate these two abelian squares.
An element of order $2$ that exchanges the first and last positions can act to obtain \emph{reap} from \emph{pear}.
On the other hand, we need an element of order $4$ to get from \emph{team} to \emph{mate}.

In Section~\ref{sec:def} we review basic concepts and notions.
We then formally define how we will use group actions to study repetitions and factor complexity.
General properties of how actions of subgroups work are discussed in Section~\ref{sec:gen}.
Connections to Strumian words are also explored in Section~\ref{sec:gen}.
Next, we study properties of $\vtm$.
In Section~\ref{sec:sqfree}, we prove strengthenings of the fact that $\vtm$ is a squarefree word.
In Section~\ref{sec:complexity}, we compute a factor complexity formula for $\vtm$, using a particular group action, in terms of the abelian complexity of the period doubling word $\pd$.

\section{Preliminaries and definitions}
\label{sec:def}
Let $\Sigma$ be a finite set. 
We call $\Sigma$ an \emph{alphabet} and call any element of $\Sigma$ a \emph{letter}. 
A \emph{word} is any sequence, finite or infinite, of letters from $\Sigma$.
We let $\Sigma^{\ell}$ denote the set of words of length $\ell$.
If $\w \in \Sigma^{\ell}$ we write $|\w| = \ell$.
We will index words from $0$ with brackets; so, if $\w \in \Sigma^{\ell}$, then $\w[i] \in \Sigma$ for each $0 \leq i \leq \ell-1$.
For any $a \in \Sigma$ we let $|\w|_a$ denote the size of the set $\{i : \w[i] = a\}$.
Also, we use $\Sigma^{\omega}$ to denote the set of infinite words $\w$ where $\w[i] \in \Sigma$ for $i \geq 0$.
A word $\u$ is a \emph{factor} of $\w$ if $\u[i] = \w[j+i]$ for some $j$ and all $0 \leq i < |\u|$; equivalently, $\u$ is a contiguous subsequence of $\w$.

A \emph{morphism} is a map $\phi$ between the set of words on two alphabets so that $\phi(\u\v) = \phi(\u)\phi(\v)$.
A morphism is determined by its output on each letter of the alphabet.
A \emph{fixed point} of a morphism $\phi$ is a word $\w$ such that $\phi(\w) = \w$.
We will make use of standard constructions of infinite words as fixed points of a morphism.
One such word is $\vtm \in \{0, 1, 2\}^{\omega}$ which is the fixed point of the nonuniform morphism defined by $0 \mapsto 012$, $1 \mapsto 02$, and $2 \mapsto 1$.
This infinite ternary word is a notable and recurring example in combinatorics on words.
A related word is the \emph{period-doubling word} $\pd \in \{0,1\}^{\omega}$, and this word is defined to be the infinite binary word that is the fixed point of the morphism given by $0 \mapsto 01$ and $1 \mapsto 00$.

Now that we have introduced some basics on words, we add a group acting on the alphabet.
Consider a group $G$ that acts on $\Sigma$, then $G$ acts on words over $\Sigma$ by acting on each letter of a word.
For example, let $\Sigma = \{0,1,2\}$, $G$ be the group of permutations of $\Sigma$, and $g \in G$ be the permutation of $\Sigma$ that is the $3$-cycle sending $i\in \Sigma$ to $i+1 \pmod 3$.
We then have $g \cdot (0102) = 1210$.
The action of the group does not affect the length of a word; therefore, the group $G$ acts on the set of words $\Sigma^{\ell}$ for any $\ell$.

Next, for each $\ell > 0$ we choose a group $H_{\ell}$ that acts on $\Sigma^{\ell}$ by permuting positions.
For example, let $\ell = 4$ and $\Sigma = \{0,1,2\}$.
Also, let $H_4$ be the group of permutations of ${0,1,2,3}$, which indexes the positions of a word of length 4.
If $h \in H_4$ is the permutation that is the product of the transposition exchanging $0$ and $1$ with the transposition exchanging $2$ and $3$, then $h \cdot (0102) = 1020$.

We then have an action of the product $G \times H_{\ell}$ on $\Sigma^{\ell}$ where
\[(g,h) \cdot \w = g\cdot(h \cdot \w) = h \cdot (g \cdot \w)\]
for each $(g,h) \in G \times H_{\ell}$ and $\w \in \Sigma^{\ell}$.
We will consider concepts in combinatorics on words using orbits of such a group action to determine equivalence classes.
In general, one must specify the groups $G$ and $H_\ell$, along with their respective actions, to define equivalence classes via orbits.

The first concept we examine is that of powers of words.
A \emph{$(G \times H_{\ell})$-$r$th power} is defined to be a word of the form $\x_1 \x_2 \cdots \x_r$ with $\x_i \in \Sigma^{\ell}$ for each $1 \leq i \leq r$, and for all for all $1 \leq i< j \leq r$ there exists $(g,h) \in G \times H_{\ell}$ so that $(g,h) \cdot \x_i = \x_j$.
If $r=2$, we use the term \emph{square} for a 2nd power.
As before, specifying the groups alone is insufficient, and we will make it clear in context how they are acting.
Usual powers come from taking $(G, H_{\ell}) = (1,1)$ while \emph{abelian powers} are realized by taking $(G, H_{\ell}) = (1, S_{\ell})$ where $S_{\ell}$ is the symmetric group permuting positions.
Let $\rot_{\ell}$ denote the cyclic group of order $\ell$ that acts on positions by cyclically permuting them.
So, the generator of $\rot_{\ell}$ is the permutation of $\{0,1,\dots, \ell-1\}$ sending $i$  $i+1$ taken modulo $\ell$ for each $0 \leq i < \ell$.
We find that a $(1, \rot_\ell)$-square is a word that is either a square or \emph{mesosome} in the sense of~\cite{mesosome}.
If one takes the binary alphabet with $G = S_2$ whose generator swaps the two letters of the alphabet and $H_{\ell} = 1$ for all $\ell > 0$, then a $(G, H_{\ell})$-square is a square or \emph{antisquare} as defined in~\cite{anti} and further studied in~\cite{myWalnut}.

The second concept we consider in our setting is that of the factor complexity of an infinite word, which counts the number of distinct factors of a given length.
For an infinite word $\w$ we let $\rho_{\w,\ell}^{(G,H_\ell)}$ denote the number of equivalence classes in $\{\u \in \Sigma^\ell : \u \subset \w\}$ under the $G \times H_{\ell}$ action.
For brevity, we use the notation $\rho_{\w,\ell}^{(G,H_\ell)}$ even though this does not explicitly include information on how the groups are acting.
When necessary, we will provide descriptions of the group actions in statements of results.
The behavior of $\rho_{\w,\ell}^{(G,H_\ell)}$ with $G = 1$ has been considered in detail by Charlier--Puzynina--Zamboni~\cite{GroupComplexity}.
The usual factor complexity of a word is the case that $G = 1$ and $H_\ell = 1$, and it will be denoted $\rho_{\w,\ell}^{(1,1)}  = \rho_{\w,\ell}$.
Taking $G = 1$ and $H_{\ell} = S_{\ell}$ acting by permuting positions, we recover the much studied \emph{abelian complexity} which we will denote $\rho_{\w,\ell}^{(1,S_{\ell})}  = \rho^{ab}_{\w,\ell}$.
Other variations of factor complexity have been considered that correspond to particular choices of $G$ and $H_{\ell}$.
A recent notion in the literature is the \emph{reflection complexity}~\cite{reflection} which corresponds to $G = 1$ and $H_{\ell}$ being the group of order $2$, which we will denote $\rev_{\ell}$, whose generator acts by reversing the order of symbols of a length $\ell$ word.
Another example is \emph{cyclic complexity}~\cite{cyclicComplexity} where $G=1$ and $H_{\ell} = \rot_{\ell}$ is a cyclic group of order $\ell$ acting by conjugation.

\section{Some general observations}
\label{sec:gen}

\subsection{Actions of subgroups}
A guiding principle is that actions of larger groups imply stronger results in terms of avoiding powers.
The following proposition makes explicit this idea, and its proof is immediate from our definition.

\begin{prop}
    Let $G' \leq G$ and $H'_{\ell} \leq H_{\ell}$ for some $\ell > 0$ be subgroups.
    If a word $\w$ avoids $(G, H_{\ell})$-$r$th powers, then $\w$ avoids $(G', H'_{\ell})$-$r$th powers
    \label{prop:subpower}
\end{prop}

Proposition~\ref{prop:subpower} is not ``strict'' in the sense that the property of avoiding powers of all lengths with larger groups acting can be equivalent to avoiding them with subgroups acting.
This can be illustrated with the aforementioned mesosomes~\cite{mesosome} as we now discuss.
For each $\ell > 0$, take $G=1$ and $H_{\ell} = \rot_{\ell}$, where we recall that the latter group acts by cyclically permuting positions (i.e., conjugation).
A $(1, \rot_{\ell})$-square is word $\x \x'$ such that $\x = \u\v$ and $\x' = \v\u$.
If $\v = \epsilon$, then $\x = \x'$ and $\x\x'$ is a square.
Otherwise, $\x\x'$ contains the square $\v^2$.
Therefore, containing a square of some length is equivalent to containing a $(1, \rot_{\ell})$-square for some $\ell > 0$.

In terms of factor complexity, we can make a similar observation.
Acting by larger groups causes our version of factor complexity to decrease.
This again follows readily from our definition.
Indeed, a larger group acting makes our orbits larger, meaning we may end up counting fewer distinct equivalence classes.

\begin{prop}
    If $G' \leq G$ and $H'_{\ell} \leq H_{\ell}$ for each $\ell > 0$ be subgroups, then
    \[\rho_{\w, \ell}^{(G, H_{\ell})} \leq \rho_{\w, \ell}^{(G', H'_{\ell})} \]
    for any $\w$ and $\ell > 0$.
    \label{prop:subcomplexity}
\end{prop}

Proposition~\ref{prop:subcomplexity} gives us an inequality that we can easily show is not strict in general.
For a simple example illustrating this, consider the word
\[(012)^{\omega} = 012012012012\cdots\]
and for $\ell \geq 2$ let $H_{\ell}$ be the group of order $2$ generated by the transposition exchanging $0$ and $1$ (i.e., exchanging the two leftmost positions of a word).
Then $\rho_{(012)^{\omega},\ell}^{(1,H_{\ell})} = \rho_{(012)^{\omega},\ell}$ for all $\ell \geq 2$.
For any factor $\u \subset (012)^{\omega}$ with $|\u| = \ell \geq 2$ the generator of $H_{\ell}$ sends $\u$ to a word of length $\ell$ that is not a factor of $(012)^{\omega}$.
The words
\begin{align*}
    (012)^{\ell} && (120)^{\ell} && (201)^{\ell}
\end{align*}
are the only three factors of $(012)^{\omega}$ with length $3\ell > 0$.
For the generator $h \in H_{3\ell}$ we have
\begin{align*}
    h((012)^{\ell}) &= (102)(012)^{\ell-1}\\
    h((120)^{\ell}) &= (210)(120)^{\ell-1}\\
    h((201)^{\ell}) &= (021)(201)^{\ell-1}
\end{align*}
none of which are factors of $(012)^{\omega}$.
Factors of other lengths can be analyzed in a similar manner.

Many more examples of equality in Proposition~\ref{prop:subcomplexity} can be found.
One can see~\cite{avoidRev} for constructions of words $\w$ satisfying $\rho_{\w,\ell}^{(1,\rev_{\ell})} = \rho_{\w,\ell}$, for sufficiently large $\ell$, where $\rev_{\ell}$ is the group of order two acting on words of length $\ell$ by reversing the order of symbols.

We can also have the maximal possible reduction in the complexity when acting by a larger group.
Let $G = S_{\{0,2\}}$ be the group of order $2$ acting on $\Sigma = \{0,1,2\}$ whose generator exchanges $0$ and $2$
In~\cite[Lemma 6]{vtm} it is shown that when $u \subset \vtm$, then also $g\cdot u \subset \vtm$ where $g$ is the generator of $S_{\{0,2\}}$.
This means that
\[\rho_{\vtm, \ell}^{(S_{\{0,2\}},1)} = \frac{\rho_{\vtm, \ell}}{2}\]
for any $\ell \geq 2$, which is the maximal decrease one can have from the standard factor complexity with a group of order $2$.
Note the statement fails for $\ell = 1$ because we have the three factors $0$, $1$, and $2$.
However, it holds for $\ell \geq 2$ because $\vtm$ is squarefree implying every factor of length $\ell \geq 2$ has either a $0$ or a $2$.

\subsection{Sturmian words}

A \emph{Sturmian word} is an infinite binary word $\w$ such that $\rho_{\w, \ell} = \ell + 1$ for all $\ell > 0$.
These are aperiodic binary words of minimal factor complexity.
The interested reader can see~\cite[Chap. 2]{Lothaire} for a more in-depth introduction to these words.
Sturmian words are also aperiodic binary words of minimal abelian complexity satisfying $\rho_{\w, \ell}^{ab} = 2$ for all $\ell > 0$.
So, it is natural to ask how Sturmian words behave with factor complexities $\rho_{\w,\ell}^{(G, H_{\ell})}$.
There are results with $G=1$ characterizing Sturmian words in terms of minimal complexity~\cite{GroupComplexity}.
We will exhibit some examples with $G\neq 1$, showing that Sturmian words do not minimize complexity among all aperiodic words for all choices of group actions.

The \emph{Parikh vector} of a word is a vector that records how many occurrences of each letter are in the word.
For a word $\w \in \Sigma$, it is the vector $(|\w|_a)_{a \in \Sigma}$.
The abelian complexity counts the number of distinct Parikh vectors of factors of a given length.
If we add an action of a group $G$ acting on our alphabet, then $G$ acts by permuting positions of Parikh vectors.
In the situation where $|\Sigma| = k$ and $G = S_k$, we can pick a preferred form for Parikh vectors.
We will use the convention of Parikh vectors sorted into weakly decreasing order.

Consider the binary alphabet $\Sigma = \{0,1\}$ and let $G = S_2$ act by exchanging the letters.
For each $\ell > 0$, we take $H_{\ell} = S_{\ell}$ acting by permuting positions.
We will find that there are aperiodic binary words with smaller complexity values than Sturmian words for this choice of group actions.
Additionally, we can find aperiodic binary words, both Sturmian and non-Sturmian, that have the same values of this complexity measure.

The \emph{Thue--Morse word} $\tm$ is the fixed point of the morphism $0 \mapsto 01$ and $1 \mapsto 10$.
Its abelian complexity is known to be
\[\rho_{\tm, \ell}^{ab} = \begin{cases} 2 & \text{if } \ell \text{ is odd}\\ 3 & \text{if } \ell \text{ is even}\end{cases}\]
for each $\ell > 0$.
Moreover, we know that Parikh vectors are $\{(m+1, m), (m, m+1)\}$ when $\ell = 2m+1$ and $\{(m+2, m), (m+1,m+1), (m, m+2)\}$ when $\ell = 2m+2$.
This can be found in~\cite[Theorem 3.3]{tmab} where the authors actually characterize all binary words with this abelian complexity.
The following proposition is now immediate.

\begin{prop}
The formula
\[\rho_{\tm, \ell}^{(S_2, S_{\ell})} = \begin{cases} 1 & \text{if } \ell \text{ is odd}\\ 2 & \text{if } \ell \text{ is even}\end{cases}\]
holds for all $\ell > 0$.
\label{prop:tm}
\end{prop}

Now let $\w_{\alpha}$ be the characteristic Sturmian word with irrational slope $0 <\alpha < \frac{1}{4}$.
This restriction on $\alpha$ makes it convenient for us to construct an example of the Sturmian word with larger complexity than Proposition~\ref{prop:tm}.
The Parikh vectors arising from factors of $\w_{\alpha}$ of length $\ell$ factors are $\left(\ell - \lfloor \ell \alpha \rfloor, \lfloor \ell \alpha \rfloor\right)$ and $\left(\ell - \lceil \ell \alpha \rceil, \lceil \ell \alpha \rceil\right)$.
For $\ell \in \{1, 2, 3, 4\}$ the Parikh vectors are
\begin{align*}
    \{(1,0), (0,1)\} && \{(2,0), (1,1)\} && \{(3,0), (2,1)\} && \{(4,0), (3,1)\}
\end{align*}
where we have used $\alpha < \frac{1}{4}$.

\begin{lem}
    If $0 < \alpha < \frac{1}{4}$, then $\ell - \lceil \ell \alpha \rceil \geq \lceil \ell \alpha \rceil$ for $\ell > 4$.
    \label{lem:1/4}
\end{lem}
\begin{proof}
    Let $\lceil \ell \alpha \rceil = \ell \alpha + \epsilon$ for $0 \leq \epsilon \leq 1$.
    Assume that $\ell - \lceil \ell \alpha \rceil < \lceil \ell \alpha \rceil$.
    This implies\
    \begin{align*}
        \ell &< 2\ell \alpha + 2 \epsilon\\
        &\leq 2\ell\alpha + 2\\
        &< \frac{\ell}{2} + 2
    \end{align*}
    which is a contradiction.
\end{proof}

\begin{prop}
If $\w_{\alpha}$ is the characteristic Sturmian word with irrational slope $0 <\alpha < \frac{1}{4}$, then
\[\rho_{\w_{\alpha}, \ell}^{(S_2, S_{\ell})} = \begin{cases} 1 & \text{if } \ell=1\\ 2 & \text{if } \ell>1\end{cases}\]
holds for all $\ell > 0$.
\label{prop:Sturmian}
\end{prop}
\begin{proof}
By Lemma~\ref{lem:1/4} and inspection for small $\ell$, we see that the Parikh vectors  $\left(\ell - \lfloor \ell \alpha \rfloor, \lfloor \ell \alpha \rfloor\right)$ and $\left(\ell - \lceil \ell \alpha \rceil, \lceil \ell \alpha \rceil\right)$ are already sorted in weakly decreasing order when $\ell > 1$.
The Proposition immediately follows as $\lfloor \ell \alpha \rfloor \neq \lceil \ell \alpha \rceil$ since $\alpha$ is irrational.
\end{proof}

Proposition~\ref{prop:tm} and Proposition~\ref{prop:Sturmian} show that Sturmian words are not characterized by having minimal complexity among aperiodic words for an arbitrary choice of group action.
We will now construct a non-Sturmian word whose complexity with the $S_2 \times S_{\ell}$ action exactly matches what we found in Proposition~\ref{prop:Sturmian}.
Before this, we establish the following lemma and make use of the Walnut theorem prover~\cite{Walnut} to give a quick mechanical proof.

\begin{lem}
    For any $\ell > 1$ and $(a,b) \in \{0,1\} \times \{0,1\}$ there exists $\u \subset \tm$ such that $\u[0] = a$ and $\u[\ell-1] = b$.
    \label{lem:startend}
\end{lem}
\begin{proof}
The four choices for $(a,b)$ can be check in Walnut with
    \begin{verbatim}
        eval lem00 "A n (n>1 => Ej (T[j]=@0 & T[j+n-1]=@0))";
        eval lem01 "A n (n>1 => Ej (T[j]=@0 & T[j+n-1]=@1))";
        eval lem10 "A n (n>1 => Ej (T[j]=@1 & T[j+n-1]=@0))";
        eval lem11 "A n (n>1 => Ej (T[j]=@1 & T[j+n-1]=@1))";
    \end{verbatim}
    which returns TRUE four times.
\end{proof}

\begin{rmk}
For readers unfamiliar with Walnut, we recommend the book~\cite{ShallitBook}.
We now briefly explain the proof of Lemma~\ref{lem:startend}.
The word $\tm$ is predefined in Walnut as $\texttt{T}$.
The command
\begin{verbatim}
            eval lem00 "A n (n>1 => Ej (T[j]=@0 & T[j+n-1]=@0))";
\end{verbatim}
tells Walnut to determine if the statement
\[\forall n, n>1 \implies \exists j, \tm[j]=0 \wedge \tm[j+n-1]=0\]
is true or false.
Here \texttt{lem00} is a user chosen name that can be used to look up output or refer to later.
\end{rmk}

\begin{prop}
If $\w = \psi(\tm)$ where $\psi$ is the morphism defined by $0 \mapsto 0011$ and $1 \mapsto 0110$, then
\[\rho_{\w, \ell}^{(S_2, S_{\ell})} = \begin{cases} 1 & \text{if } \ell=1\\ 2 & \text{if } \ell>1\end{cases}\]
holds for all $\ell > 0$. Futhermore, $\w$ is not a Sturmian word.
\label{prop:nonSturmian}
\end{prop}
\begin{proof}
Our word $\w = \psi(\tm)$ is constructed by concatenations of $\psi(0) = 0011$ and $\psi(1) = 0110$, which we will refer to as blocks.
Each block has a Parikh vector of $(2,2)$.
It follows $\big||\u|_0 - |\u|_1\big|$ is bounded that for any factor $\u \subset \w$.
This absolute difference can be analyzed by considering $|\u| \pmod 4$ along with looking at prefixes and suffixes of our blocks $\psi(0)$ and $\psi(1)$.

Let us consider a factor $\u \subset \w$ expressed as $\u = \A\B\C$ where $\B$ is a concatenation of some number of blocks while $\A$ and $\C$ are a suffix and a prefix of a block, respectively.
Here $\A$ and $\C$ may be empty, but neither is a complete block because in that case it would be part of $\B$.
First assume $|\u| = 4m$.
All possible cases are $\u = \B$ or $\u$ is one of the following
\begin{align*}
    1\B001 && 1\B011 && 0\B001 && 0\B011\\
    011\B0 && 011\B0 && 110\B0 && 110\B0
\end{align*}
implying the possible Parikh vectors of $\u$ are $(2m, 2m)$, $(2m+2, 2m-2)$, and $(2m-2, 2m+2)$.

Now consider $|\u| = 4m + 1$, and $\u$ must be one of the following
\begin{align*}
    0\B && 1\B && \B0 && \B0\\
    11\B001 && 11\B011 && 10\B001 && 10\B011\\
    011\B00 && 011\B01 && 110\B00 && 110\B01\\
\end{align*}
resulting in $(2m+1, 2m)$, $(2m, 2m+1)$, $(2m+2, 2m-1)$, and $(2m-1, 2m+2)$ as possible Parikh vectors.

Similar analysis shows that when $|\u| = 4m+2$ its Parikh vector is among $(2m+1, 2m+1)$, $(2m+2, 2m)$, and $(2m,2m+2)$.
Lastly, when $|\u| = 4m+3$ the possible Parikh vectors are $(2m+2, 2m+1)$, $(2m+1, 2m+3)$, $(2m+3,2m)$ and $(2m,2m+3)$.
Thus, we find that the formula in the proposition is an upper bound for $\rho_{\w, \ell}^{(S_2, S_{\ell})}$.
The equality follows from Lemma~\ref{lem:startend}, which tells us that all possible prefixes $\A$ and suffixes $\C$ from our case checking occur.
To see that $\w$ is not Sturmian we observe that its abelian complexity function is not constant and equal to $2$.
\end{proof}

\section{Strengthening squarefreeness of $\vtm$}\label{sec:sqfree}

One of the key properties of $\vtm$ is that it is a squarefree word, which means $\vtm$ avoids $(1,1)$-squares in our notation for any length $\ell$.
In this section, we look to strengthen the idea of $\vtm$ being squarefree by showing instances where $\vtm$ avoids squares using our new definition with nontrivial groups acting.

If we have $|\Sigma| = k$ and $G = \mathbb{Z}/k\mathbb{Z}$ acting by cyclically permuting $\Sigma$ we will call $(\mathbb{Z}/k\mathbb{Z}, 1)$-$r$th power a \emph{Caesar $r$th power} where here we take $H_{\ell} = 1$ for any $\ell > 0$.
We use this terminology in reference to the classic \emph{Caesar shift cipher} as a Caesar $r$th power is a word $\x_1 \x_2 \dots \x_r$ such that each $\x_i$ is a possible encryption of $\x_1$ with such a cipher.
With a binary alphabet of a Caesar $2$nd power, or \emph{Caesar square}, is a binary word that is either a square or an antisquare.
In this way, we see that the cyclic action of the alphabet is one possible extension of the idea of binary complementation to larger alphabets.

If a group $G$ acts transitively on $\Sigma$, then any length $r$ word over $\Sigma$ is a $(G \times H_1)$-$r$th power.
So, in the case of $G$ acting transitively, it is impossible to construct an infinite word avoiding $(G \times H_{\ell})$-$r$th powers for all $\ell$ and fixed $r$.
Note that in the case of Caesar squares, we have $G = \mathbb{Z}/k\mathbb{Z}$ acting transitively on our alphabet.
So, as is sometimes done when studying repetition avoidance, we will avoid all sufficiently long Caesar squares.
Another approach was taken in~\cite{NewProb} where the authors instead opted not to use the entire group, but instead avoided factors of the form 
\[xx, xg(x), xg^2(x), \dots, xg^j(x)\]
for some chosen $j$ where $g$ is a generator of the group.

Now consider ternary words with $G =\mathbb{Z}/3\mathbb{Z}$.
Explicitly, we have the group generator acting on the alphabet $\{0,1,2\}$ by sending $i$ to $i+1 \pmod{3}$ for $0 \leq i < 3$.
We will focus on finding squarefree ternary words such that they avoid Caesar squares of sufficiently large length.
By directly checking, we can find that every sufficiently long ternary squarefree word contains some Caesar square of length less than or equal to $10$.
Up to permutation of our ternary alphabet, the unique longest ternary squarefree word avoiding Caesar squares of length $10$ or greater is
\[010201210212012102010\]
which has length $21$.
It turns out that $\vtm$ is a ternary squarefree word avoiding Caesar squares of length $12$ or greater.
As the discussion in this paragraph shows, this result is optimal in terms of finding infinite ternary squarefree words avoiding ``short'' Caesar squares.
We again make use of Walnut~\cite{Walnut, ShallitBook} to prove this result.

\begin{thm}
The squarefree ternary word $\vtm$ avoids Caesar squares of length $12$ or greater.
\label{thm:Caesar}
\end{thm}
\begin{proof}
Running the commands
\begin{small}
\begin{verbatim}
eval no_sq0 "~ Ej En Ai (n > 0) & ((i < n) => (VTM[j+i] = VTM[j+i+n]))";

eval no_sq1 "~ Ej En Ai (n > 5) & ((i < n) => ((VTM[j+i]=@0 & VTM[j+i+n]=@1) | 
    (VTM[j+i]=@1 & VTM[j+i+n]=@2) | (VTM[j+i]=@2 & VTM[j+i+n]=@0)))";

eval no_sq2 "~ Ej En Ai (n > 5) & ((i < n) => ((VTM[j+i]=@0 & VTM[j+i+n]=@2) | 
    (VTM[j+i]=@1 & VTM[j+i+n]=@0) | (VTM[j+i]=@2 & VTM[j+i+n]=@1)))";
\end{verbatim}
\end{small}
in Walnut we see \texttt{TRUE} three times, confirming the theorem is true.
The \texttt{eval no\_sq0} checks the well-known fact that $\vtm$ is squarefree.
The \texttt{eval no\_sq1} and \texttt{eval no\_sq2} check that $\vtm$ does not contain a Caesar square of length $12$ or greater with a group element acting send $0$ to $1$ or $0$ to $2$ respectively.
\end{proof}

We may wonder if $\vtm$ avoids $(S_3,1)$-squares of lengths greater than $N$ for some value of $N$.
Here $S_3$ is the symmetric group of all permutations of our alphabet $\{0,1,2\}$.
It turns out that $\vtm$ contains arbitrarily long $(S_3,1)$-squares.

\begin{figure}
    \centering
    \includegraphics[width=0.9\linewidth]{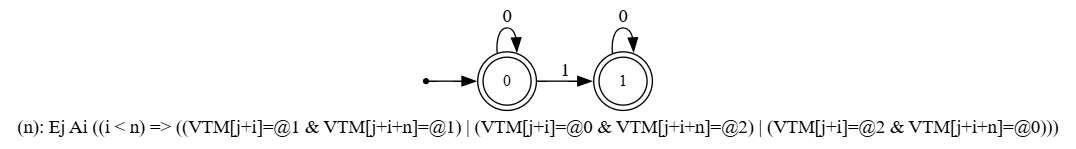}
    \caption{An automaton accepting binary strings that represent lengths of certain $(S_3,1)$-squares in $\vtm$.}
    \label{fig:AUT}
\end{figure}

\begin{prop}
   The squarefree ternary word $\vtm$ contains arbitrarily long $(S_3,1)$-squares.
\end{prop}
\begin{proof}
    We can prove this statement again using Walnut.
    If we type
    \begin{verbatim}
        eval long_S_sq "Ej Ai ((i < n)=>((VTM[j+i]=@1 & VTM[j+i+n]=@1) 
                                | (VTM[j+i]=@0 & VTM[j+i+n]=@2) 
                                | (VTM[j+i]=@2 & VTM[j+i+n]=@0)))";
    \end{verbatim}
    we obtain the automaton depicted in Figure~\ref{fig:AUT}.
    In the Walnut command, the variable $n$ is unquantified.
    This means that Walnut will give us an automaton that accepts all values of $n$, in base-$2$, making the statement true.
    Here $n$ represents half the length of a $(S_3,1)$-square that makes use of the permutation of $\{0,1,2\}$ that swaps $0$ and $2$ while fixing $1$.
    We can see that this automaton accepts all powers of $2$, and therefore $\vtm$ contains arbitrarily long $(S_3,1)$-squares.
\end{proof}

In Theorem~\ref{thm:Caesar} we saw $\vtm$ avoided all Caesar squares of sufficiently long length.
Now, we wish to focus on cases where we do not need to specify the length of our avoidance phenomenon.
Since $\vtm$ contains any ternary word of length $2$ as a factor, we cannot avoid $(G, H_1)$-squares if $G$ is nontrivial.
So, we investigate if $\vtm$ can avoid $(1, H_{\ell})$-squares for all $\ell > 0$ where the group $H_{\ell}$ is not trivial.

Recall, $\rev_{\ell} \cong \Z/2\Z$ is the group whose generator acts on length $\ell$ words by reversing the order of the symbols.
So, the generator of $\rev_5$ would act on $11001$ by sending it to $10011$.
We say that a \emph{reflection square} is a word $w = xy$ so that $x = g \cdot y$ where $|x| = \ell = |y|$ and $g$ is the generator of $\rev_{\ell}$.
This means that a $(1, \rev_{\ell})$-square is either a square or a reflection square.

\begin{thm}
    For any $\ell > 0$, $\vtm$ avoids $(1, \rev_{\ell})$-squares. That is, $\vtm$ contains neither reflection squares nor squares.
\end{thm}
\begin{proof}
    We may appeal to Walnut and run
    \begin{verbatim}
        eval no_sq_rev "~ Ej En Ai (n > 0) 
                &  ((i < n) => (VTM[i+j] = VTM[i+2*n-j-1]))";
    \end{verbatim}
    which returns \texttt{TRUE}. Thus, the theorem is proven as we already know that $\vtm$ avoids squares.
\end{proof}

One can find further examples of this kind of square avoidance, like the following theorem.

\begin{thm}
    For any $\ell > 0$ let $H_{\ell} = \Z/2\Z$ act on words of length $\ell$ by its generator exchanging the leftmost and rightmost letters of the word, then $\vtm$ avoids $(1, H_{\ell})$-squares.
\end{thm}

\begin{proof}
We already know that $\vtm$ avoids usual squares.
So, it suffices to run
    \begin{verbatim}
        eval no_sq_first_last "~ Ej En Ai (n > 0) 
            & (VTM[i] = VTM[i+2*n-1]) & (VTM[i+n-1] = VTM[i+n]) 
            & ((i > 0) & (i < n-1) => (VTM[j+1] = VTM[j+i+n]))";
    \end{verbatim}
    in Walnut that verifies we have no $(1, H_{\ell})$-squares with a nontrivial action of the group for any $\ell > 0$.
\end{proof}

\section{Factor Complexity of $\vtm$ with Action of Symmetric Groups}\label{sec:complexity}

We now wish to analyze the quantity $\rho_{\vtm,\ell}^{(S_3,S_\ell)}$ which we will relate to the abelian complexity $\rho_{\pd,\ell}^{ab}$.
Here $S_3$ is a symmetric group permuting our ternary alphabet, and $S_{\ell}$ is a symmetric group permuting the positions of a length $\ell$ factor of $\vtm$.
It will be advantageous to review some relevant facts about the words $\vtm$ and $\pd$, which can be found in~\cite[Section 3]{vtm}.
Firstly, if we remove all $1$'s from $\vtm$ we obtain the infinite word $(02)^{\omega}  = 020202\cdots$ and thus for any factor $\w \subset \vtm$ satisfies 
\begin{equation}
\Big| |\w|_0 - |\w|_2 \Big| \leq 1
\label{eq:bal}
\end{equation}
giving us a balance property on the number of $0$'s relative to the number of $2$'s in any factor of $\vtm$.
If $\phi$ is the morphism defined by $\phi(0) = 0 = \phi(2)$ and $\phi(1) = 1$, then $\phi(\vtm) = \pd$.
For a given $\ell > 0$ the set $\{|\w|_1 : \w \subset \vtm, |\w| = \ell\}$ is an interval of consecutive integers without any gaps.

It will be important to understand the quantities
\[p_{\ell} = \min \{|\w|_1 : \w \subset \vtm, |\w| = \ell\}\]
and
\[q_{\ell} = \max \{|\w|_1 : \w \subset \vtm, |\w| = \ell\}\]
giving the minimum and maximum number of $1$'s that occur in a length $\ell$ factor of $\vtm$.
We have the recurrences
\begin{align}
    p_{4\ell - 1} &= p_{\ell} + \ell - 1\label{eqpa}\\
    p_{4\ell} &= p_{\ell} + \ell\\
    p_{4\ell + 1} &= p_{\ell} + \ell\\
    p_{4\ell + 2} &= p_{\ell} + \ell\label{eqpb}
\end{align}
and
\begin{align}
    q_{4\ell - 1} &= q_{\ell} + \ell \label{eqa}\\
    q_{4\ell} &= q_{\ell} + \ell\\
    q_{4\ell + 1} &= q_{\ell} + \ell + 1\\
    q_{4\ell + 2} &= q_{\ell} + \ell + 1\label{eqb}
\end{align}
which were written down in this form in~\cite{KSZ}. The initial conditions are $p_1 = p_2 = 0$ and $q_1 = q_2 = 1$.

\begin{lem}
    For any $\ell > 0$ the inequality 
    \[p_{\ell} \leq \left\lfloor \frac{\ell - 1}{3} \right\rfloor\]
    holds.
    \label{lem:floor3p}
\end{lem}

\begin{proof}
    The lemma follows by induction using the recurrences for $p_{\ell}$ found in Equations~(\ref{eqpa})-(\ref{eqpb}).
    Let us assume that
\[p_{\ell} \leq \left\lfloor \frac{\ell - 1}{3} \right\rfloor\]
while also observing that we know this to be true for $\ell=1$ and $\ell = 2$.
    We then have
    \begin{align*}
        p_{4\ell - 1} &= p_{\ell} + \ell - 1\\
        &\leq \left\lfloor \frac{\ell - 1}{3} \right\rfloor + \ell - 1\\
        &= \left\lfloor \frac{4\ell - 4}{3} \right\rfloor\\
        &\leq \left\lfloor \frac{4\ell - 2}{3} \right\rfloor
\end{align*}
and the inequality holds when $\ell \pmod{4} \equiv -1$.
    We also find that
\begin{align*}
        p_{4\ell+2 } &= p_{4\ell + 1}\\
        &= p_{4\ell} \\
        &= p_{\ell} + \ell\\
        &\leq \left\lfloor \frac{\ell - 1}{3} \right\rfloor + \ell\\
        &= \left\lfloor \frac{4\ell - 1}{3} \right\rfloor\\
        &\leq \left\lfloor \frac{4\ell}{3} \right\rfloor\\
        &\leq \left\lfloor \frac{4\ell+1}{3} \right\rfloor
\end{align*}
and the inequality holds when $\ell \pmod{4} \not\equiv -1$.
\end{proof}

\begin{lem}
    For any $\ell > 0$ the inequality 
    \[ q_{\ell} \geq \left\lfloor \frac{\ell}{3} \right\rfloor + 1\]
    holds.
    \label{lem:floor3q}
\end{lem}
\begin{proof}
    The lemma follows by induction using the recurrences for $q_{\ell}$ found in Equations~(\ref{eqa})-(\ref{eqb}).
    Let us assume that
    \[q_{\ell} \geq \left\lfloor \frac{\ell}{3} \right\rfloor + 1\]
    while also observing we know this to be true for $\ell=1$ and $\ell = 2$.
    We then have
    \begin{align*}
        \left\lfloor \frac{4\ell-1}{3}\right\rfloor + 1&\leq \left\lfloor \frac{4\ell}{3}\right\rfloor + 1\\
        &= \left\lfloor \frac{3\ell + \ell}{3} \right\rfloor + 1\\
        &= \left\lfloor \frac{\ell}{3} \right\rfloor + 1 + \ell\\
        &\leq q_{\ell} + \ell\\
        &= q_{4\ell}\\
        &= q_{4\ell-1}
    \end{align*}
    and the inequality is satisfied when $\ell \pmod{4}$ is $-1$ or $0$.
    We also have
    \begin{align*}
        \left\lfloor \frac{4\ell+1}{3}\right\rfloor + 1 &\leq \left\lfloor \frac{4\ell+2}{3}\right\rfloor + 1\\
        &= \left\lfloor \frac{3\ell + \ell + 2}{3} \right\rfloor + 1\\
        &= \left\lfloor \frac{\ell + 2}{3} \right\rfloor + \ell + 1\\
        &\leq \left\lfloor \frac{\ell}{3} \right\rfloor + 1 + \ell + 1\\
        &\leq q_{\ell} + \ell + 1\\
        &= q_{4\ell + 2}\\
        &= q_{4\ell + 1}
    \end{align*}
   and the inequality is satisfied when $\ell \pmod{4}$ is $2$ or $3$.
\end{proof}

We record the following lemma for later use.
Its proof is straightforward.

\begin{lem}
    Let $a$, $b$, $c$, and $d$ be integers such that $a+b = c + d$.
    If $|a-b| \leq 1$ and $|c-d| \leq 1$, then $\{a,b\} = \{c, d\}$
    \label{lem:abcd}
\end{lem}

We now prove the main theorem of this section.

\begin{thm}
   The equality
   \[\rho_{\vtm,\ell}^{(S_3,S_\ell)} = \rho_{\pd,\ell}^{ab} - 1\]
   holds for any $\ell > 0$.
\end{thm}

\begin{proof}
For each orbit of factors of $\vtm$ under the $S_3 \times S_{\ell}$, we will analyze which abelian equivalence class or classes of $\pd$ the orbit maps to under $\phi$.
Take nonnegative integers $a$, $b$, and $c$ with $a+b+c = \ell$.
Assume there exists a factor $\u \subset \vtm$ with Parikh vector $(|\u|_0, |\u|_1, |\u|_2) = (a,b,c)$.
This implies that $|a-c| \leq 1$ by Equation~(\ref{eq:bal}).

First, assume that $|a-b| > 1$ and $|b-c| > 1$.
This means that the $S_3$ action on the alphabet $\{0,1,2\}$ cannot exchange $1$ with either of $0$ or $2$ while staying within the set of factors of $\vtm$.
In this case, an equivalence class of the factors of $\vtm$ under the $S_3 \times S_{\ell}$ action directly corresponds with an abelian equivalence class of factors of $\pd$.
Here $\phi(\u) \subset \pd$ has $(|\phi(\u)|_0 , |\phi(\u)|_1) = (a+c,b)$.
Any $\u' \subset \vtm$ in the orbit of $\u$ under the $S_3 \times S_{\ell}$ action will also have $(|\phi(\u')|_0 , |\phi(\u')|_1) = (a+c,b)$.
Furthermore, if $\v \subset \vtm$ with $(|\phi(\v)|_0 , |\phi(\v)|_1) = (a+c,b)$, then $\v$ must be in the orbit of the $S_3 \times S_{\ell}$ as $\u$. This follows from Lemma~\ref{lem:abcd} since $|\v|_0 + |\v|_2 = |\u|_0 + |\u|_2$ with $\big||\v|_0 - |\v|_2\big| \leq 1$ and $\big||\u|_0 - |\u|_2\big| \leq 1$.

Next, assume that $|a-b| \leq 1$ or $|b-c| \leq 1$.
In this case, we will show there is an orbit of factors of $\vtm$ under the $S_3 \times S_{\ell}$ action which maps to two different abelian equivalence classes of factors of $\pd$.
If $\ell \equiv 1 \pmod 3$ the only options for $(a,b,c)$ are $(a+1,a,a)$, $(a,a+1,a)$ and $(a,a,a+1)$.
Consider $\u,\u',\u'' \subset \vtm$ with
\begin{align*}
    (|\u|_0, |\u|_1, |\u|_2) &= (a+1,a,a)\\
    (|\u'|_0, |\u'|_1, |\u'|_2) &= (a,a+1,a)\\
    (|\u''|_0, |\u''|_1, |\u''|_2) &= (a,a,a+1)\\
\end{align*}
so that $\u$, $\u'$, and $\u''$ are all equivalent under the $S_3 \times S_{\ell}$ action.
However, $\phi(\u')$ is in a different abelian equivalence class of factors on $\pd$ than $\phi(\u)$ and $\phi(\u'')$.
The case that $\ell \equiv 2 \pmod{3}$ is similar.

The case that $\ell \equiv 0 \pmod{3}$ has factors with Parikh vectors $(a,a,a)$ in one $S_3 \times S_{\ell}$ orbit of factors of $\vtm$ mapping to its own abelian equivalence class of $\pd$.
We then also have Parikh vectors $(a,a+1,a-1)$ and $(a-1, a+1, a)$ along with $(a,a-1,a+1)$ and $(a+1,a-1,a)$ all coming from one $S_3 \times S_{\ell}$ orbit of factors of $\vtm$ while mapping to two abelian equivalence classes of $\pd$ with Parikh vectors $(2a-1,a+1)$ and $(2a+1,a-1)$ respectively.

Since we know $\phi(\vtm) = \pd$ we can conclude that $\rho_{\vtm,\ell}^{(S_3,S_\ell)} \leq \rho_{\pd,\ell}^{ab}$.
To conclude that $\rho_{\vtm,\ell}^{(S_3,S_\ell)} = \rho_{\pd,\ell}^{ab} - 1$ we need only know that the case in the previous paragraph with a single $S_3 \times S_{\ell}$ orbit of factors of $\vtm$ mapping to two distinct abelian equivalence classes of $\pd$ occurs for any $\ell$.
By Lemma~\ref{lem:floor3p} and Lemma~\ref{lem:floor3q} this does occur and the theorem is proven.
\end{proof}

Now that we know explicitly the close relationship between $\rho_{\vtm,\ell}^{(S_3,S_\ell)}$ and $\rho_{\pd,\ell}^{ab}$ we can easily deduce further facts about $\rho_{\vtm,\ell}^{(S_3,S_\ell)}$ since the abelian complexity $\rho_{\pd,\ell}^{ab}$ is well understood.
For example, it is known that
\begin{align*}
    \rho_{\pd,2\ell}^{ab} &= \rho_{\pd,\ell}^{ab}\\
    \rho_{\pd,4\ell \pm 1}^{ab} &= \rho_{\pd,\ell}^{ab} + 1
\end{align*}
from~\cite[Proposition 2]{vtm}. It follows that
\begin{align*}
    \rho_{\vtm,2\ell}^{(S_3,S_{2\ell})} &= \rho_{\vtm,\ell}^{(S_3,S_\ell)}\\
    \rho_{\vtm,4\ell \pm 1}^{(S_3,S_{4\ell \pm 1})} & = \rho_{\vtm,\ell}^{(S_3,S_\ell)} + 1
\end{align*}
must also hold.
This means that in the Online Encyclopedia of Integer Sequences the values of $\rho_{\vtm,\ell}^{(S_3,S_\ell)}$ are~\cite[\href{https://oeis.org/A007302}{A007302}]{OEIS}.

\bibliographystyle{alphaurl}
\bibliography{refs}

\end {document}